\documentclass[12pt]{article}

\usepackage{latexsym}
\usepackage{amssymb}
\usepackage{amsmath}
\usepackage{amsthm}
\usepackage{amscd}
\usepackage[pdftex]{graphicx}
\usepackage{comment}

\newtheorem{thm}{Theorem}[section]
\newtheorem{cor}[thm]{Corollary} 
\newtheorem{lem}[thm]{Lemma} 
\newtheorem{prop}[thm]{Proposition} 
\newtheorem{rem}[thm]{Remark} 
\newtheorem{exam}[thm]{Example} 
 
\newtheorem{defin}[thm]{Definition}

\newcommand{\enD}{\hfill $\Box$\vspace{3truemm} \par}
\newcommand{\R}{\mathbb{R}}

\newcommand{\bmu}{\mbox{\boldmath $\mu$}}

\newcommand{\bx}{\mbox{\boldmath $x$}}

\newcommand{\A}{\mathcal{A}}

\def\proof{\noindent {\sl Proof} :\;  }

\begin{document}

\title{Singularities in bivariate normal mixtures}

\author{Yutaro Kabata, Hirotaka Matsumoto, Seiichi Uchida and Masao Ueki}

\date{\today}

\maketitle
\begin{abstract}
We investigate mappings \( F = (f_1, f_2) \colon \mathbb{R}^2 \to \mathbb{R}^2 \) where \( f_1, f_2 \) are bivariate normal densities from the perspective of singularity theory of mappings, motivated by the need to understand properties of two-component bivariate normal mixtures. We show a classification of mappings \( F = (f_1, f_2) \) via $\mathcal{A}$-equivalence and characterize them using statistical notions. Our analysis reveals three distinct types, each with specific geometric properties. Furthermore, we determine the upper bounds for the number of modes in the mixture for each type.
\end{abstract}

\renewcommand{\thefootnote}{\fnsymbol{footnote}}

\section{Introduction}
A normal mixture is an important statistical model frequently used to represent multimodal distributions in real-world data analysis. This paper aims to study two-component bivariate normal mixtures from the perspective of the singularity theory of mappings.

\begin{figure}[ht]
\begin{center}
\includegraphics[clip,width=10.0cm]{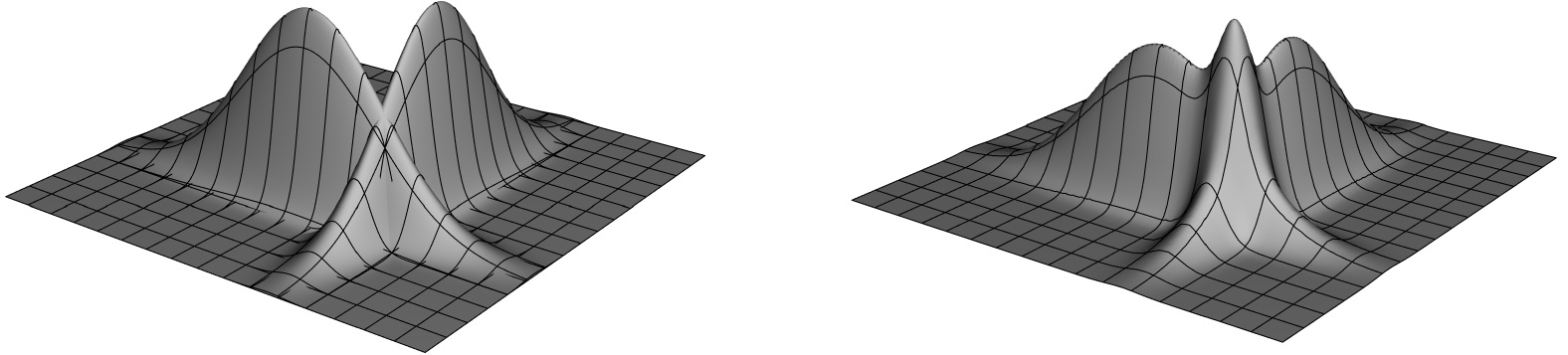}
\end{center}
\caption{The left figure is the graph of two normal densities $f_1, f_2$ with $\bmu_1=(0,0)$, $\bmu_2=(1,1)$, $\Sigma_1=\left(\begin{array}{cc}1&0\\0&0.2\\\end{array}\right)$, $\Sigma_2=\left(\begin{array}{cc}0.2&0\\0&1\\\end{array}\right)$. The right figure is the graph of the mixture density $\frac12 f_1+\frac12 f_2$ which has three modes (local maximum points).}
\label{fig:exceeding mixture}
\end{figure}

Let \( f_i \colon \mathbb{R}^2 \to \mathbb{R}_{>0} \) for \( i = 1, 2 \) denote the densities of bivariate normal distributions, defined as
\begin{eqnarray}\label{eq:normaldis}
f_i(\bx) := \phi(\bx; \boldsymbol{\mu}_i, \Sigma_i) = \frac{1}{2\pi |\Sigma_i|^{\frac{1}{2}}} \exp \left( -\frac{1}{2} (\bx - \boldsymbol{\mu}_i)^{T} \Sigma_i^{-1} (\bx - \boldsymbol{\mu}_i) \right),
    \end{eqnarray}
where \(\bx = (x, y)\) is a variable, \(\boldsymbol{\mu}_i \in \mathbb{R}^2\) is a mean vector, and \(\Sigma_i\) is a positive definite \(2 \times 2\) covariance matrix. Given two bivariate normal densities \(f_1\) and \(f_2\), the density of normal mixture \(M_c\) is expressed as their convex linear combination:
\[
M_c := c f_1 + (1 - c) f_2,
\]
where \(0 \leq c \leq 1\). Despite \(M_c\) being a simple linear combination of two functions, its behavior is highly nontrivial. In particular, the number of modes (local maximum points) of \(M_c\) is of significant interest in statistics and varies depending on the parameters \(c\), \(\boldsymbol{\mu}_i\), and \(\Sigma_i\). The maximum possible number of modes for the mixture, known as its \textit{modality}, has been studied extensively in the literature \cite{Amendola2020, CarreiraPerpinan2003, RayRen2012, RayLindsay2005, Lindsay1983}. 
 Figure \ref{fig:exceeding mixture} illustrates a typical example where a
 two-component bivariate normal mixture can exhibit three modes, exceeding the number of components. This example serves as a somewhat surprising counterexample that contradicts the intuition of practitioners in fields such as statistics, machine learning and image processing. {Therefore, it is both theoretically intriguing and practically significant to inquire about the conditions under which the number of modes in a mixture distribution does not exceed the number of components. Building on this, this study provides a novel contribution to this problem by examining mixture distributions from the perspective of singularity theory of mappings.}

The mixture density can be expressed as
\[
M_c = c f_1 + (1 - c) f_2 = (c, 1 - c) \cdot (f_1, f_2),
\]
where \( M_c \) is the inner product of the vector \( (c, 1 - c) \in \mathbb{R}^2 \) and the mapping \( F = (f_1, f_2) \colon \mathbb{R}^2 \to \mathbb{R}_{>0}^2 \). This relationship suggests that studying the mapping \( F \) provides significant insights into the properties of \( M_c \).
In particular, properties such as modality, which do not depend on the value of the mixing proportion \( c \), can be obtained by studying \( F \). Indeed, it has been implied that the inflection points on the boundary curve of the image of \( F \) are closely related to the number of modes of \( M_c \) \cite{RayLindsay2005}. 
It should be noted that the above boundary curve can be regarded as singularities of \( F \).

Based on singularity theory of mappings (cf. \cite{IRFT}), the present paper investigates the geometry of the product mapping \( F = (f_1, f_2) \colon \mathbb{R}^2 \to \mathbb{R}_{>0}^2 \) of normal distributions \( f_1 \) and \( f_2 \). 
We denote {\it the singular set of $F$} by $S(F):=\{\bx\in\R^2\;|\; \det J_F(\bx)=0 \}$, and {\it the singular value set of $F$} by $C(F):=F\left(S(F)\right)$. Note that the boundary curve of the image of a mapping $F$ is contained in $C(F)$. Since topological properties of $S(F)$ and $C(F)$ do not change under the coordinate changes, we consider the classification of mappings via the following {\it $\A$-equivalence}. Two smooth mappings $F,G\colon\R^2\to\R^2$ are said to be $\A$-equivalent if and only if there exist diffeomorphisms $\Phi\colon\R^2\to\R^2$ of the source space and $\Psi\colon\R^2\to\R^2$ of the target space so that $\Psi\circ F\circ \Phi =G$. 

Table \ref{tab:classification} summarizes our results.
We show that, via the $\A$-equivalence, the space of product mappings $F=(f_1,f_2)$ are divided into only three $\A$-equivalent classes under a natural assumption $\bmu_1\not=\bmu_2$. Each class is clearly characterized by the singular set $S(F)$ and types of singularities (Corollary \ref{cor:classificationofF}). 

\begin{table}[!htb]
\begin{center}
\begin{tabular}{c|c|c|c|c|c}
\hline
Type & $S(F)$&Singularities&Proportional&Codirectional&Modalities\\
&&&&\\
\hline
\hline
1& hyperbola & $(x,y^2)$,&no&no&1, 2 \mbox{or} 3\\
&&$(x,xy+y^3)$&&\\
\hline
2& two intersecting lines & $(x,y^2)$, &no&yes &1 \mbox{or} 2\\
&&$(x,xy^2+y^4)$&&\\
\hline
3& line & $(x,y^2)$&yes&$-$& 1 \mbox{or} 2\\
&&&&\\
\hline
\end{tabular}
\end{center}
\caption{Classification of the pairs of normal densities $f_1,f_2$ ($\bmu_1\not=\bmu_2$). The third column shows the local normal forms of possible singularities of $F=(f_1,f_2)$. The first row of the fourth column indicates whether the covariance matrices of the two density functions are proportional, while the second row indicates whether they are codirectional. The fifth column means the possible modalities of the mixture distributions $M_c=cf_1+(1-c)f_2$ for each type.}
\label{tab:classification}
\end{table}

The above classification is given in terms of singularity theory of a mapping. On the other hand, the classification can be characterized in rather statistical or linear algebraic notions for the pair of normal densities: {\it proportionality of covariance matrices} and {\it codirectionality} (see \S 3). 
Also, the upper bounds of the number of modes of the mixture $M_c$ for each type are determined (see \S 4). 
{Types 2 and 3 represent cases where the number of modes does not exceed the number of components. These types are distinctly characterized by the proportionality of covariance matrices and codirectionality. In particular, Type 2 has, to the best of our knowledge, not been previously addressed in the literature.}

The paper is organized as follows. In Section 2, we introduce the concept of {\it a generalized distance-squared mapping} and its classification via \( A \)-equivalence, following \cite{Ichiki2016}. In Section 3, we provide a classification of the product mappings of two bivariate normal distributions, based on the results of Section 2. Additionally, the concept of codirectionality is introduced to characterize the above classification. In Section 4, the number of the modes of the mixture is investigated for each type in the classification.

\begin{rem}{\rm
The concept of using mappings to analyze mixtures has its origins in earlier works such as \cite{RayLindsay2005, Lindsay1983}. This paper builds upon these foundational ideas and extends them by applying perspectives from singularity theory of mappings to further investigate these approaches. It is worth noting that an additional advantage of this research direction is the ability to easily visualize properties of mixture distributions by drawing the image of the mapping, as demonstrated in \S5. This visualization technique provides a powerful tool for understanding and analyzing complex mixture models.

Beyond the study of mixtures in statistics, the relationship between a scalar function and a mapping is also a key topic in the field of multi-objective optimization problems. Notably, approaches from differential topology and singularity theory, as explored in \cite{hamada2020, smale1973, smale1974, smale1975}, offer valuable insights into this subject.
}
\end{rem}

\section{Generalized distance-squared mapping}
Let ${\bf p}_1 = (p_{11}, p_{12}), {\bf p}_2 = (p_{21}, p_{22})\in\R^2$, \( A = (a_{ij})_{1 \leq i \leq 2, 1 \leq j \leq 2} \) be a \(2 \times 2\) matrix with non-zero entries. Then the following mapping \( G({\bf p}_1, {\bf p}_2, A) : \mathbb{R}^{2} \to \mathbb{R}^{2} \) is called a generalized distance-squared mapping:
  \begin{eqnarray*}
G({\bf p}_1, {\bf p}_2, A)(x,y) := \left( a_{11}(x-p_{11})^2+a_{12}(y-p_{12})^2, a_{21}(x-p_{21})^2+a_{22}(y-p_{22})^2
 \right).
  \end{eqnarray*}
The generalized distance-squared mapping is introduced and investigated in terms of $\A$-equivalence in \cite{Ichiki2016}.
 In \S 3, we show that the mapping \( F = (f_1, f_2) \colon \mathbb{R}^2 \to \mathbb{R}^2 \) with bivariate normal distributions \( f_1, f_2 \) is $\A$-equivalent to a generalized distance-squared mapping. Thus we quickly review the results of \cite{Ichiki2016} in this section. Note that there is an exceptional case not addressed in \cite{Ichiki2016}, which we cover in Proposition \ref{prop:generalizedSQexception}.
  
We denote \( A_k \) by a \( 2 \times 2 \) matrix of rank \( k \) with non-zero entries. We sum up the results of \cite{Ichiki2016} used in the present paper as follows.

\begin{prop}[{\cite{Ichiki2016}}]\label{prop:normalform}
Suppose ${\bf p}_{1}\not={\bf p}_{2}$. Then the following hold:
\begin{enumerate}
\item The mapping \( G({\bf p}_1, {\bf p}_2, A_1) \) is \(\A\)-equivalent to \(\ (x,y^2)\).
\item The mapping \( G({\bf p}_1, {\bf p}_2, A_2) \) is \(\A\)-equivalent to
$$
\left((x - q)^2 + (y - r)^2, ax^2 + by^2\right),
$$
where $(q,r)\not= (0,0)$ and $a,b>0$ hold. 
\end{enumerate}
\end{prop}
\proof
Giving suitable coordinate changes are enough.
For the statement 1, see the ``Proof of part (1) of Theorem 1" in \cite{Ichiki2016}. For the statement 2, see the ``Proof of Proposition 3" in \cite{Ichiki2016}.
 Note that, although the statements ``(1) of Theorem 1" and ``Proposition 3" in \cite{Ichiki2016} is stated under the condition with respect to $({\bf p}_1,{\bf p}_2)$: $p_{11}\not=p_{21}$ and $p_{12}\not=p_{22}$, the coordinate changes given in their proofs only need the condition ${\bf p}_1\not={\bf p}_2$. Thus the statements and their proofs of Proposition \ref{prop:normalform} make sense for general $({\bf p}_1,{\bf p}_2)$ with the coordinate changes given in \cite{Ichiki2016}.

\begin{thm}[Theorem 1 and Proposition 2 in \cite{Ichiki2016}]
 Suppose $p_{11}\not=p_{21}$ and $p_{12}\not=p_{22}$. Then the following hold:
  \begin{enumerate}
    \item The mapping \( G({\bf p}_1, {\bf p}_2, A_1) \) is \(\A\)-equivalent to \(\ (x,y^2)\).
    \item The singular set \( S(G({\bf p}_1, {\bf p}_2, A_2)) \) is a rectangular hyperbola. Any point of $S(G({\bf p}_1, {\bf p}_2, A_2))$ is a fold point except for one; and the exceptional point is a cusp. In particular, 
for any \(2 \times 2\) matrix $\tilde{A}_2$ with non-zero entries and rank $2$, \( G({\bf p}_1, {\bf p}_2, A_2) \) is $\A$-equivalent to \( G({\bf p}_1, {\bf p}_2, \tilde{A}_2) \).
     \end{enumerate}
     \label{thm:generalizedSQ}
\end{thm}

Note that in the above statement, a point $\bx_0\in\R^2$ is called a fold point (resp. cusp) of the mapping $F\colon\R^2\to\R^2$ if $F$ is locally expressed as $(x,y^2)$ (resp. $(x,xy+y^3)$) by taking local coordinate changes of the source space $\R^2$ around $\bx_0$ and the target space around $F(\bx_0)$ (for the detail of singularity theory, see \cite{IRFT} for example).

For our purpose, we need  the following statements which deal with the exceptional cases of Theorem \ref{thm:generalizedSQ}.

\begin{prop}\label{prop:generalizedSQexception}
Suppose that either $p_{11}=p_{21}$ or $p_{12}=p_{22}$, but not both, holds. Then the following hold:
  \begin{enumerate}
    \item The mapping \( G({\bf p}_1, {\bf p}_2, A_1) \) is \(\A\)-equivalent to \(\ (x,y^2)\).
    \item The mapping \( G({\bf p}_1, {\bf p}_2, A_2) \) is \(\A\)-equivalent to $(x,xy^2+y^4)$. The singular set \( S(G({\bf p}_1, {\bf p}_2, A_2)) \) is  two intersecting lines. Any point of \( S(G({\bf p}_1, {\bf p}_2, A_2)) \) is a fold point except for one (the node); and the singular value set \( C(G({\bf p}_1, {\bf p}_2, A_2)) \)  is a union of a smooth curve and a double point curve. 
     \end{enumerate}
\end{prop}
\proof
The first statement  immediately follows from the statement 1 of Proposition \ref{prop:normalform}.

We show the second statement. From the statement 2 of Proposition \ref{prop:normalform} and the assumption with respect to ${\bf p}_1,{\bf p}_2$, \( G({\bf p}_1, {\bf p}_2, A_2) \) is \(\A\)-equivalent to
$$
H_1(x, y) = ((x-q)^2 + y^2, ax^2 + by^2),
$$
where $q, a,b$ are nonzero constants and $a\not=b$. By the diffeomorphism of the target space $(X,Y)\mapsto (X-\frac ba Y-\frac ba q^2,Y)$, $H_1$ is $\A$-equivalent to
$$
H_2(x,y)=\left(\left(1-\frac ab\right)(x-q)^2, ax^2+by^2\right).
$$
By routine coordinate changes of the source and target spaces, we have the following equivalences:
\begin{eqnarray*}
H_2(x,y)&\sim_\A& (x^2, a(x+q)^2+by^2)\\
&\sim_\A& (x^2, 2aqx+by^2)\\
&\sim_\A& \left(x^2, x+\frac{b}{2aq}y^2\right)=:H_3(x,y)\\
\end{eqnarray*}
By replacing $x$ by $x-\frac{b}{2aq}y^2$ and then taking routine coordinate changes, we have the following:
\begin{eqnarray*}
H_3(x,y)&\sim_\A& \left(x^2-\frac{b}{aq}xy^2+\frac{b^2}{4a^2q^2}y^4, x\right)\\
&\sim_\A&\left(x, -\frac{b}{aq}xy^2+\frac{b^2}{4a^2q^2}y^4\right)\\
&\sim_\A&\left(x, xy^2+y^4\right)=:U(x,y).
\end{eqnarray*}
It is easily checked that the singular set of the map $U$ is
$$
S(U)=\{y=0\}\cup\{x+2y^2=0\}
$$
and the singular value set is the union of (i) a line parameterized by $(x,0)$ and (ii) a double point curve parameterized by $(-2y^2,y^4)$.
\enD

\begin{rem}{\rm
The map germ $f\colon \R^2,0\to\R^2,0; (x,y)\mapsto(x,xy^2+y^4)$ is finitely $\mathcal{K}$-determined but not finitely $\mathcal{A}$-determined (see Lemma 3.2.1:1 in \cite{Rieger1987}). Thus the above mapping $U(x,y)=(x, xy^2+y^4)$ has a very degenerate singularity at the origin. In addition, \cite{Farnik2017} shows the classification of mappings $F=(f_1,f_2)\colon\R^2\to\R^2$ with $f_1,f_2$ being quadratic polynomials via affine transformations of the source and target spaces. The type denoted by ``$f_3=(x^2+y,y^2)$" in \cite{Farnik2017} corresponds to our mapping $U(x,y)=(x, xy^2+y^4)$.}
\end{rem}


\section{Product mapping of bivariate normal distributions}
We set $f_i\colon \R^2\to\R$ for $i=1,2$ as the densities of bivariate normal distributions written as (\ref{eq:normaldis}). We are interested in $\A$-equivalent class of the product mapping $F:=(f_1,f_2)\colon\R^2\to\R_{>0}^2$.

\begin{cor}
Suppose $\bmu_1\not=\bmu_2$.
Then there are only three $\A$-equivalent classes for $F=(f_1,f_2)\colon\R^2\to\R_{>0}^2$, and
one of the following holds:
\begin{enumerate}
\item  $S(F)$ is a hyperbola. Any point of $S(F)$ is a fold point except for one; and the exceptional point is a cusp. {The pair of densities of this case corresponds to Type 1 in Table \ref{tab:classification}.}
\item  $S(F)$ is two intersecting lines. Any point of $S(F)$ is a fold point except for one (the node); and the singular value set $C(F)$ is a union of a smooth curve and a double point curve. In particular, $F$ is \(\A\)-equivalent to $(x,xy^2+y^4)$. {The pair of densities of this case corresponds to Type 2 in Table \ref{tab:classification}.}
\item  $S(F)$ is a line. Any point of $S(F)$ is a fold point. In particular, $F$ is \(\A\)-equivalent to \(\ (x,y^2)\). {The pair of densities of this case corresponds to Type 3 in Table \ref{tab:classification}.}
\end{enumerate}
\label{cor:classificationofF}
\end{cor}
\proof 
It is enough to consider the $\A$-equivalent class of $F$. First, $F$ is $\A$-equivalent to the following mapping $\tilde{F}$ whose components are positive definite quadratics:
\begin{equation}
\tilde{F}(x,y)=\left((\bx-\bmu_1)^T\Sigma_1^{-1}(\bx-\bmu_1),
(\bx-\bmu_2)^T\Sigma_2^{-1}(\bx-\bmu_2)
\right),
\end{equation}
which is given by the coordinate change $(X,Y)\mapsto(\log X,\log Y)$ on the target space of $F$. Furthermore, by suitable affine transformations of the source space of $\tilde{F}$, $\tilde{F}$ is $\A$-equivalent to a generalized distance-squared mapping: First, take the affine transformation so that $\tilde{F}_1(x,y)=(x-\mu_{11})^2+(y-\mu_{12})^2$; second, take the affine transformation (the composition of a rotation around $(\mu_{11},\mu_{12})$ and a translation) so that $\tilde{F}_2(x,y)=ax^2+by^2$ for nonzero constants $a, b$.

Then, according to Theorem \ref{thm:generalizedSQ} and Proposition \ref{prop:generalizedSQexception}, we have the statement.
 \enD
\begin{rem}{\rm
In our setting with just two components and two variables, the product mappings of the densities of normal distributions or positive definite quadratic forms are $\A$-equivalent to generalized distance-squared mappings. However, this does not happen in general. For example, $H=(h_1,h_2,h_3)\colon\R^2\to\R_{\ge0}^3$ with $h_1,h_2,h_3$ being bivariate normal distributions is not $\A$-equivalent to a generalized distance-squared mapping in general.}
\end{rem}

The above Corollary \ref{cor:classificationofF} gives a classification of $F=(f_1,f_2)$ with respect to singularities. 
Each type in the above classification is characterized by statistical or linear algebraic notions: proportionality of covariance matrices and codirectionality.
Here, two matrices $\Sigma_1, \Sigma_2$ are said to be proportional if there exists a constant $c>0$ so that $\Sigma_1=c\Sigma_2$. Furthermore, the notion of codirectionality is introduced as follows:

{
\begin{defin}\label{def:codirectional}
Let $f_i(\bx)=\phi(\bx;\bmu_i, \Sigma_i)$ be densities of bivariate normal distributions for $i=1,2$ where $\Sigma_1$ and $\Sigma_2$ are not proportional. 
We say that {$f_1$ and $f_2$ are codirectional} if the vector $\bmu_1-\bmu_2$ is the eigenvector of both $\Sigma_1$ and $\Sigma_2$.
\end{defin}
}


\begin{figure}[htbp]
\begin{center}
\includegraphics[clip,width=5.0cm]{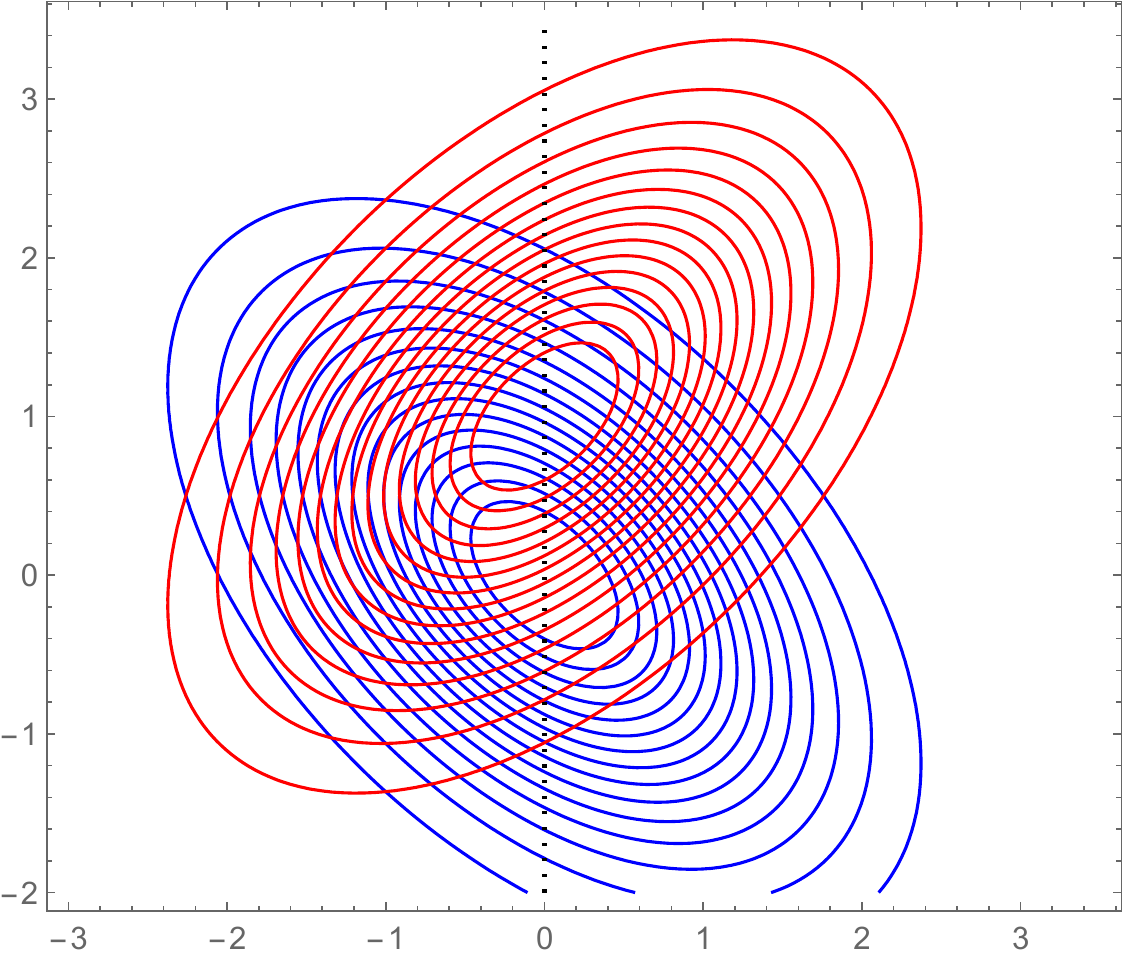}
\;\;\;
\includegraphics[clip,width=5.0cm]{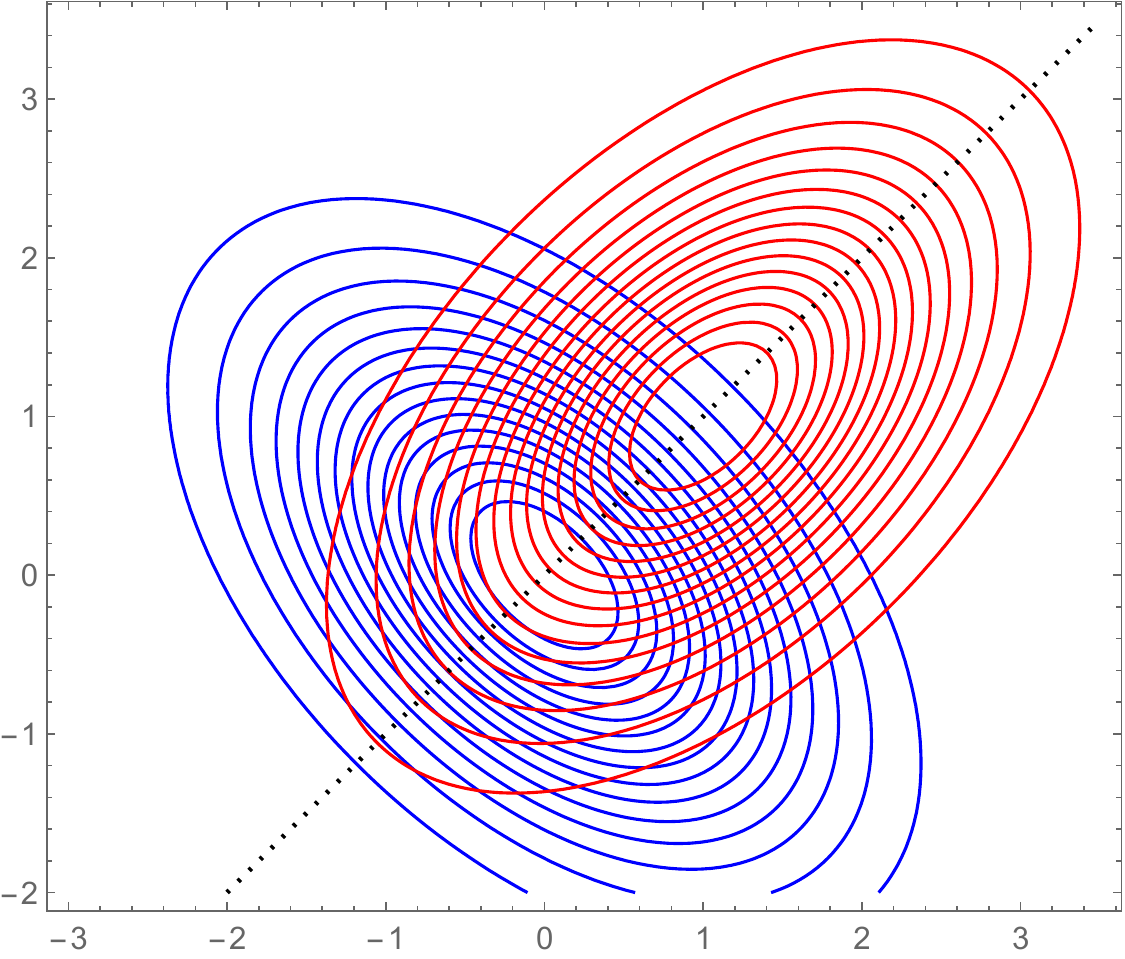}
\end{center}
\caption{The contours of the densities $f_1$ (blue) and $f_2$ (red). The dot lines show the lines going through the mean vectors $\bmu_1$ and $\bmu_2$.
The left figure shows a non-codirectional case, and the right figure shows a codirectional case.}
\label{fig:codirectionality}
\end{figure}

\begin{exam}
\begin{enumerate}
\item
Let $\bmu_1=(0,0)$, 
$\Sigma_1=\left(\begin{array}{cc}1&-0.5\\-0.5&1\\\end{array}\right)$,
$\bmu_2=(0,1)$, 
$\Sigma_2=\left(\begin{array}{cc}1&0.5\\ 0.5&1\\\end{array}\right)$.
Then $f_1, f_2$ are not codirectional.
See the left figure of Figure \ref{fig:codirectionality}.
\item
Let $\bmu_1=(0,0)$, 
$\Sigma_1=\left(\begin{array}{cc}1&-0.5\\-0.5&1\\\end{array}\right)$,
$\bmu_2=(1,1)$, 
$\Sigma_2=\left(\begin{array}{cc}1&0.5\\ 0.5&1\\\end{array}\right)$.
Then $f_1, f_2$ are codirectional.
See the right figure of Figure \ref{fig:codirectionality}.
\end{enumerate}
\end{exam}

\begin{lem}\label{lem:codirectioninvariance}
Let $f_i(\bx)=\phi(\bx;\bmu_i, \Sigma_i)$ be densities of bivariate normal distributions for $i=1,2$ where $\Sigma_1$ and $\Sigma_2$ are not proportional. Suppose $f_1$ and $f_2$ are codirectional. Then, for an affine transformation $p\colon\R^2\to\R^2$, $f_1\circ p$ and $f_2\circ p$ are codirectional.
\end{lem}
\proof
The proof immediately follows from the definition.
\enD

\begin{prop}
Set $\bmu_1=(0,0)$, $\bmu_2=(m_1,m_2)$, $\Sigma_1=I$ and $\Sigma_2=\Sigma:=\left(
\begin{array}{cc}
\sigma_1^2 &0\\
0&\sigma_2^2
\end{array}
\right)$, where $\sigma_1, \sigma_2>0$.
Then we have the following:
\begin{enumerate}
\item If $\sigma_1\not=\sigma_2$ and $m_1m_2\not=0$, then $S(F)$ is a hyperbola.
\item If $\sigma_1\not=\sigma_2$ and $m_1m_2=0$, then $S(F)$ is two intersecting lines. 
\item If $\sigma_1=\sigma_2$, then $S(F)$ is a line.
\end{enumerate} 
\label{prop:sfinnormalform}
\end{prop}
\proof
The proof follows from Theorem \ref{thm:generalizedSQ} and Proposition \ref{prop:generalizedSQexception}, since $F=(f_1,f_2)$ is regarded as a generalized distance-squared mapping by the coordinate change $(X,Y)\mapsto(\log X,\log Y)$ on the target space of $F$, in this case.
We note that the singular set $S(F)$ is defined by the following quadratic equation:
\begin{equation}
\lambda(x,y):=(x,y)Q
\left(
\begin{array}{c}
x\\
y
\end{array}
\right)
+L\left(
\begin{array}{c}
x\\
y
\end{array}
\right)=0,
\end{equation}
where
\begin{eqnarray}
Q:&=&
\left(
\begin{array}{cc}
0&-\frac12(\sigma_1^2-\sigma_2^2)\\
-\frac12(\sigma_1^2-\sigma_2^2)&0
\end{array}
\right),\\
L:&=&(m_2\sigma_1^2, -m_1\sigma_2^2).
\end{eqnarray}
In particular, $\lambda(x,y)=0$ is never an ellipse or parabola.
\enD

\begin{thm}\label{thm:characterization}
Let $f_i(\bx)=\phi(\bx;\bmu_i, \Sigma_i)$ be densities of bivariate normal distributions for $i=1,2$, and set $F=(f_1,f_2)$. The following hold:
\begin{enumerate}
\item
Suppose $\Sigma_1$ and $\Sigma_2$ are not proportional, and {$f_1$ and $f_2$ are not codirectional}. Then  $S(F)$ is a rectangular hyperbola {(Type 1)};
\item
Suppose $\Sigma_1$ and $\Sigma_2$ are not proportional, and {$f_1$ and $f_2$ are codirectional}. Then $S(F)$ is two intersecting lines {(Type 2)};
\item
Suppose $\Sigma_1$ and $\Sigma_2$ are proportional. Then $S(F)$ is a line {(Type 3)}.
\end{enumerate}
\end{thm}
\proof
It is easily seen that $F$ is $\A$-equivalent to $(\phi(\bx; 0, I),\phi(\bx; \bmu, \Sigma))$ with $\bmu=(m_1,m_2)$ and $\Sigma:=\left(
\begin{array}{cc}
\sigma_1^2 &0\\
0&\sigma_2^2
\end{array}
\right)$, where $\sigma_1, \sigma_2>0$.
In particular, $\sigma_1\not=\sigma_2$ and $m_1m_2\not=0$ hold in case 1 of Theorem \ref{thm:characterization}; $\sigma_1\not=\sigma_2$ and $m_1m_2=0$ in case 2; and $\sigma_1=\sigma_2$ in case 3.

Thus according to Lemma \ref{lem:codirectioninvariance} and Proposition \ref{prop:sfinnormalform} \ref{cor:classificationofF}, we prove the statements.
\enD

Based on the Corollary \ref{cor:classificationofF} and Theorem \ref{thm:characterization}, we have the characterization of each $\A$-equivalent class of the product mapping $F=(f_1,f_2)$ with respect to the proportionality of covariance matrices and the codirectionality as in Table \ref{tab:classification}.

\section{Modality}
In this section, we discuss the modality of the mixture $M_c=cf_1+(1-c)f_2$ for densities $f_1,f_2$ of bivariate normal distributions. The modality of two-component normal mixtures is studied in detail in \cite{RayRen2012, RayLindsay2005}. In particular, the following results are presented in these works: 

\begin{thm}[Theorem 2 in \cite{RayRen2012}]
The number of modes of $M_c$ is less than or equal to $3$.
\end{thm}

\begin{thm}[Cororllary 4 in \cite{RayLindsay2005}]\label{thm:equalsigma}
If $\Sigma_1$ and $\Sigma_2$ are proportional, then the number of modes of $M_c$ is less than or equal to $2$.
\end{thm}

\begin{rem}{\rm
In fact, results are given in more general settings of the dimension of variables and the number of components in \cite{RayRen2012, RayLindsay2005}.
}
\end{rem}

Note that the example given in Figure \ref{fig:exceeding mixture} shows that for $f_1,f_2$ of Type 1, the mixture $M_c$ can have three modes. Theorem \ref{thm:equalsigma} shows that the upper bound of the number of modes for Type 3 is two. It is natural to ask the upper bound of Type 2. To answer this question, we quickly review notions and results given in \cite{RayRen2012, RayLindsay2005}, which provides useful tools to analyze the number of modes of the mixture.

The {\it ridgeline} ${\bf x}^*\colon (0,1) \to \R^2$ is defined as
\begin{equation*}
{\bf x}^*(\alpha) = S_\alpha^{-1} \left[ (1-\alpha) \Sigma_1^{-1} \bmu_1 + \alpha \Sigma_2^{-1} \bmu_2 \right],
\end{equation*}
where $S_\alpha = \left[ (1-\alpha) \Sigma_1^{-1} + \alpha \Sigma_2^{-1} \right]$.
The ridgeline ${\bf x}^*$ is contained in $S(F)$, and any critical points of $M_c=cf_1+(1-c)f_2$ for $c\in [0,1]$ lies on it. We call the image of the ridge line by $F=(f_1,f_2)$ as {\it the image ridgeline}, and denote it by $F(x^*)$. The number of inflection points of the image ridge line where the sign of the curvature changes is crucial to the upper bounds of the modes of the mixture, and the number is equal to the zeros of the polynomial 
$$q(\alpha)= 1 - \alpha(1 - \alpha)p(\alpha),$$ 
where
\begin{equation*}
p(\alpha) = (\bmu_2 - \bmu_1)^T \Sigma_1^{-1} S_\alpha^{-1} \Sigma_2^{-1} S_\alpha^{-1}\Sigma_2^{-1}S_\alpha^{-1}\Sigma_1^{-1} (\bmu_2 - \bmu_1).
\end{equation*}

Summing up the results and discussions in \cite[Section 5]{RayLindsay2005}, we have the following claim (see also \cite[Result 1]{RayRen2012}).
\begin{thm}[\cite{RayRen2012, RayLindsay2005}]\label{thm:modesnumbercriteirion}
If $q(\alpha)$ has $n$ roots within the range $\alpha\in [0,1]$,
then the number of modes of $M_c$ is less than or equal to $\frac n2 +1$.
\end{thm}


Using the above results, we get the following Theorem \ref{thm:modestype2}, which shows the modality for Type 2.
\begin{thm}\label{thm:modestype2}
Suppose $\Sigma_1$ and $\Sigma_2$ are not proportional, and  $f_1$ and $f_2$ are codirectional. Then the number of modes of $M_c$ is less than or equal to $2$.
\end{thm}
\proof
Since the number of modes of $M_c$ are invariant under an affine transformation of the domain $\R^2$, we may assume that  $\bmu_1=(0,0)$, $\bmu_2=(m_1,0)$, $\Sigma_1=I$ and $\Sigma_2=\Sigma:=\left(
\begin{array}{cc}
\sigma_1^2 &0\\
0&\sigma_2^2
\end{array}
\right)$, where $\sigma_1> \sigma_2>0$.
In this case, we have
$$
p(\alpha)=\frac{m_1^2\sigma_1^2}{(\alpha+\sigma_1^2-\alpha\sigma_1^2)^3},
$$
and $q(\alpha)=0$ is equivalent to the following equation:
$$
(\sigma_1^2-1)^3\alpha^3-\sigma_1^2\left(m_1^2+3\left(\sigma_1^2-1\right)^2\right)\alpha^2+\sigma_1^2\left(m_1^2+3\sigma_1^2\left(\sigma_1^2-1\right)^2\right)\alpha-\sigma_1^6=0.
$$
Thus $q(\alpha)=0$ has at most three distinct solutions.
According to Theorem \ref{thm:modesnumbercriteirion}, we have the statement.
\enD

\section{Examples}
In this section, we present several examples of the contour plots of density functions $f_1,f_2$ and the images of their corresponding mappings $F=(f_1,f_2)\colon\R^2\to\R^2$ for different parameter values. Through these examples, we visually demonstrate the shapes and properties of each pair of density functions \(f_1, f_2\) and the corresponding mapping \(F\) for each type.

\subsection{Type 1}
\begin{figure}[h]
\begin{center}
\includegraphics[clip,width=6.0cm]{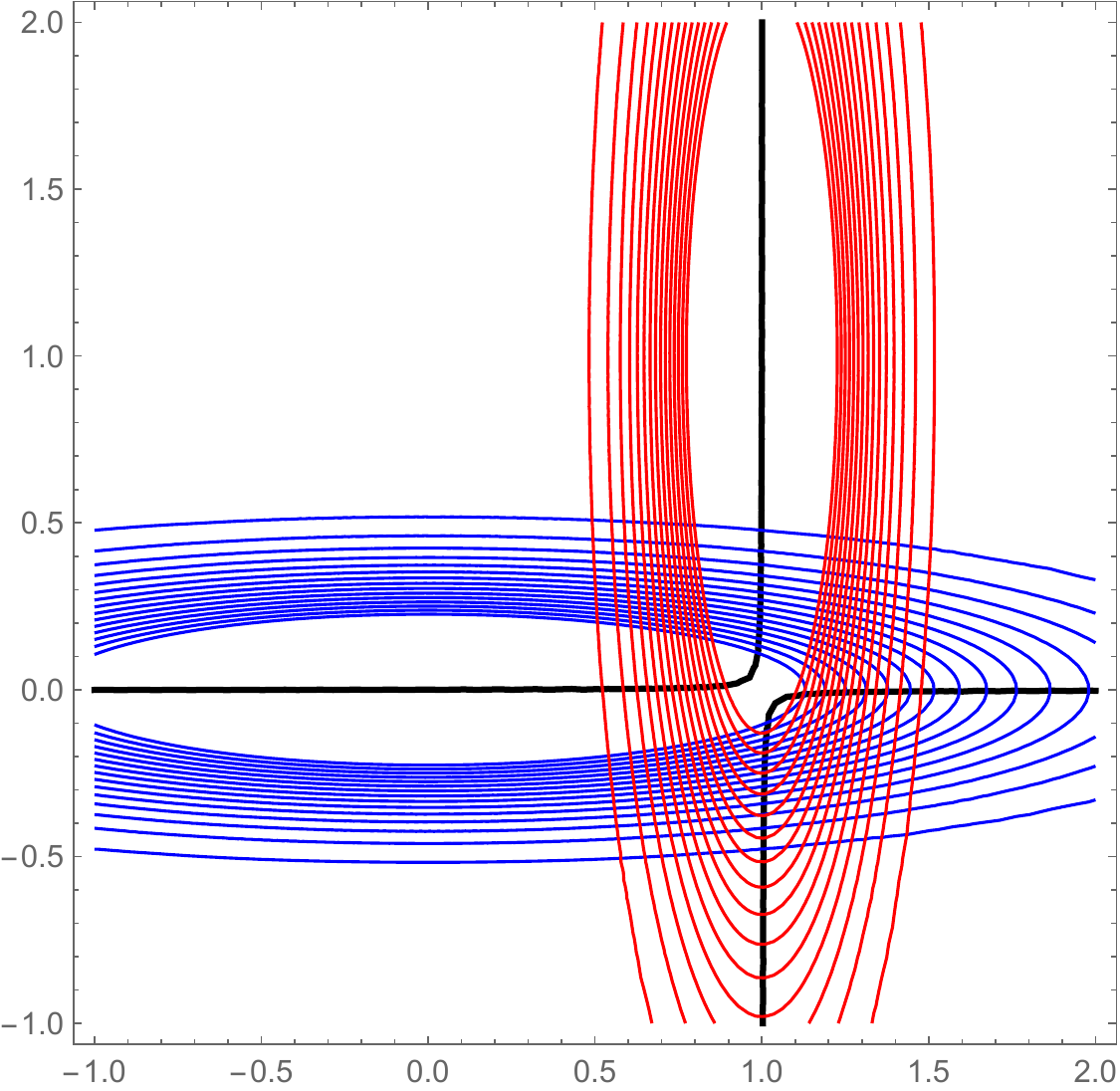}\;\;
\includegraphics[clip,width=5.9cm]{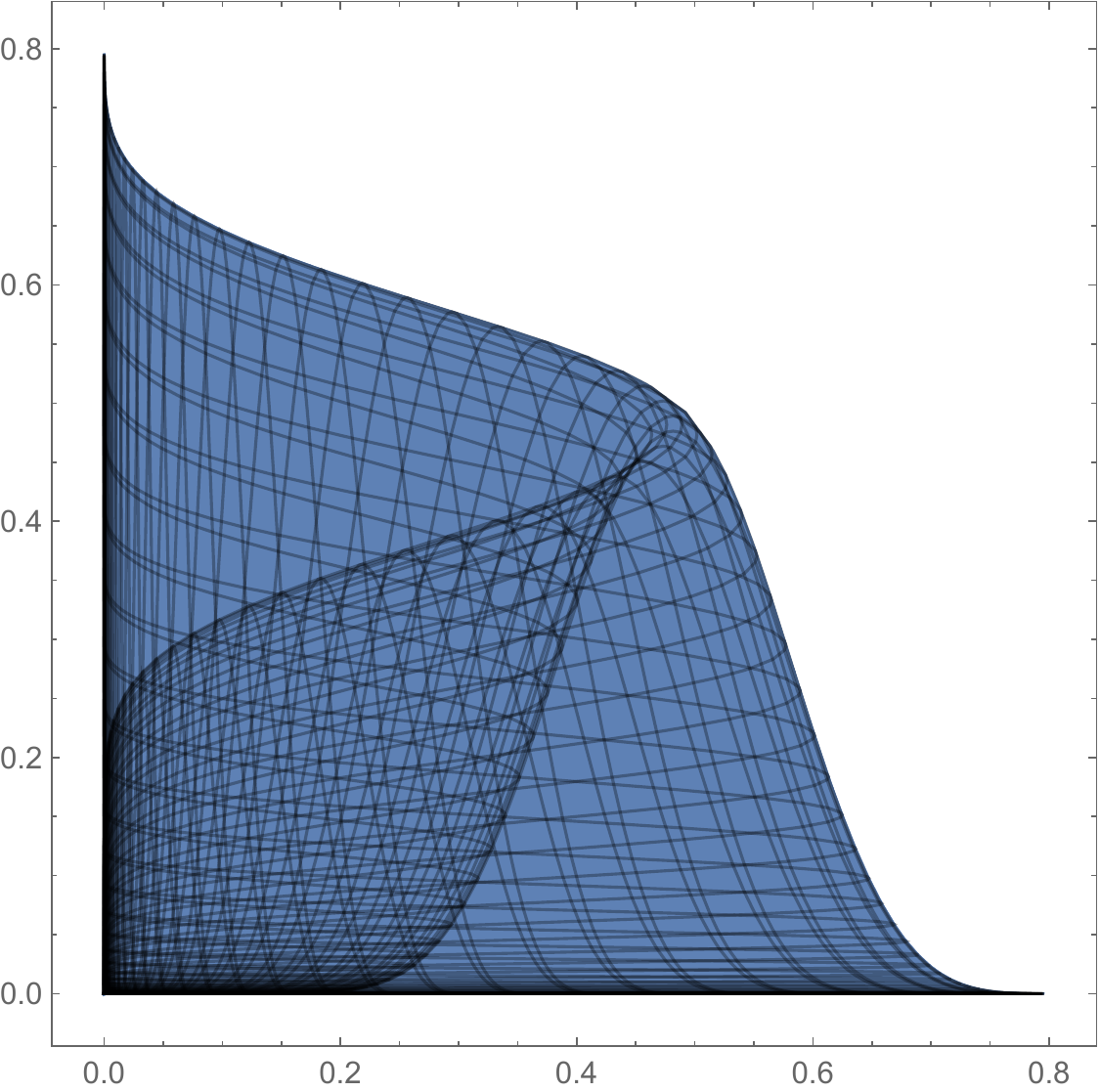}
\end{center}
\caption{Type 1: The left figure shows the contours of the densities $f_1$ (blue) and $f_2$ (red) which are not codirectional. Here the black curves are the singular set of $F=(f_1,f_2)\colon\R^2\to\R^2$. The right figure shows the image of the mapping $F\colon\R^2\to\R^2$.}
\label{fig:type1}
\end{figure}

\begin{exam}
Let
$\bmu_1=(0,0)$, 
$\Sigma_1=\left(\begin{array}{cc}1&0\\0&0.2\\\end{array}\right)$,
$\bmu_2=(1,0)$, 
$\Sigma_2=\left(\begin{array}{cc}0.2&0\\0&1\\\end{array}\right)$.
Then the pair $(f_1, f_2)$ is of Type 1. Figure \ref{fig:type1} shows the contours of $f_1,f_2$ and the image of the mapping $F\colon\R^2\to\R^2$. In particular, $S(F)$ is a hyperbola containing a unique cusp.
\end{exam}

\subsection{Type 2}
\begin{figure}[htbp]
\begin{center}
\includegraphics[clip,width=6.0cm]{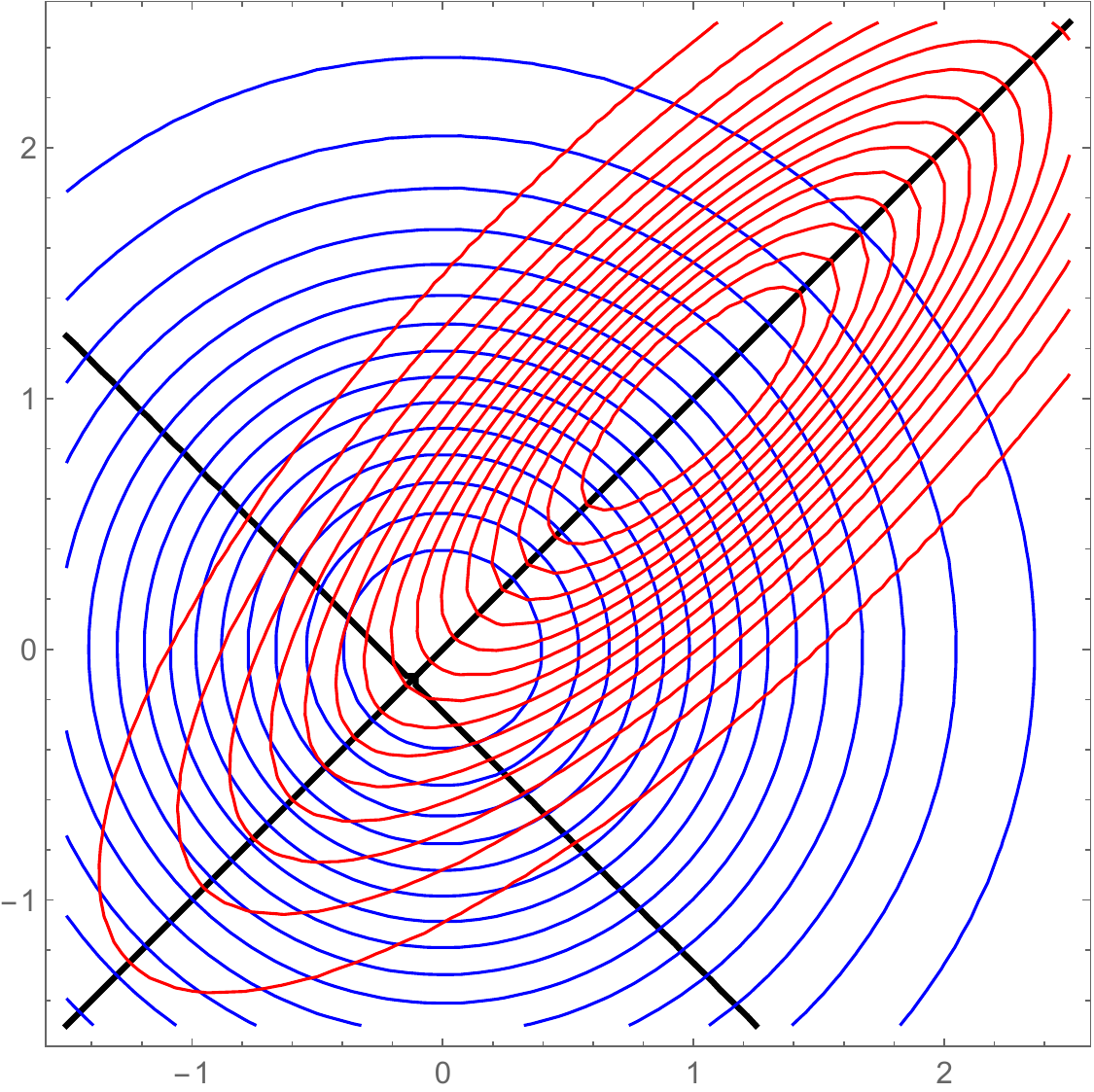}\;\;
\includegraphics[clip,width=3.9cm]{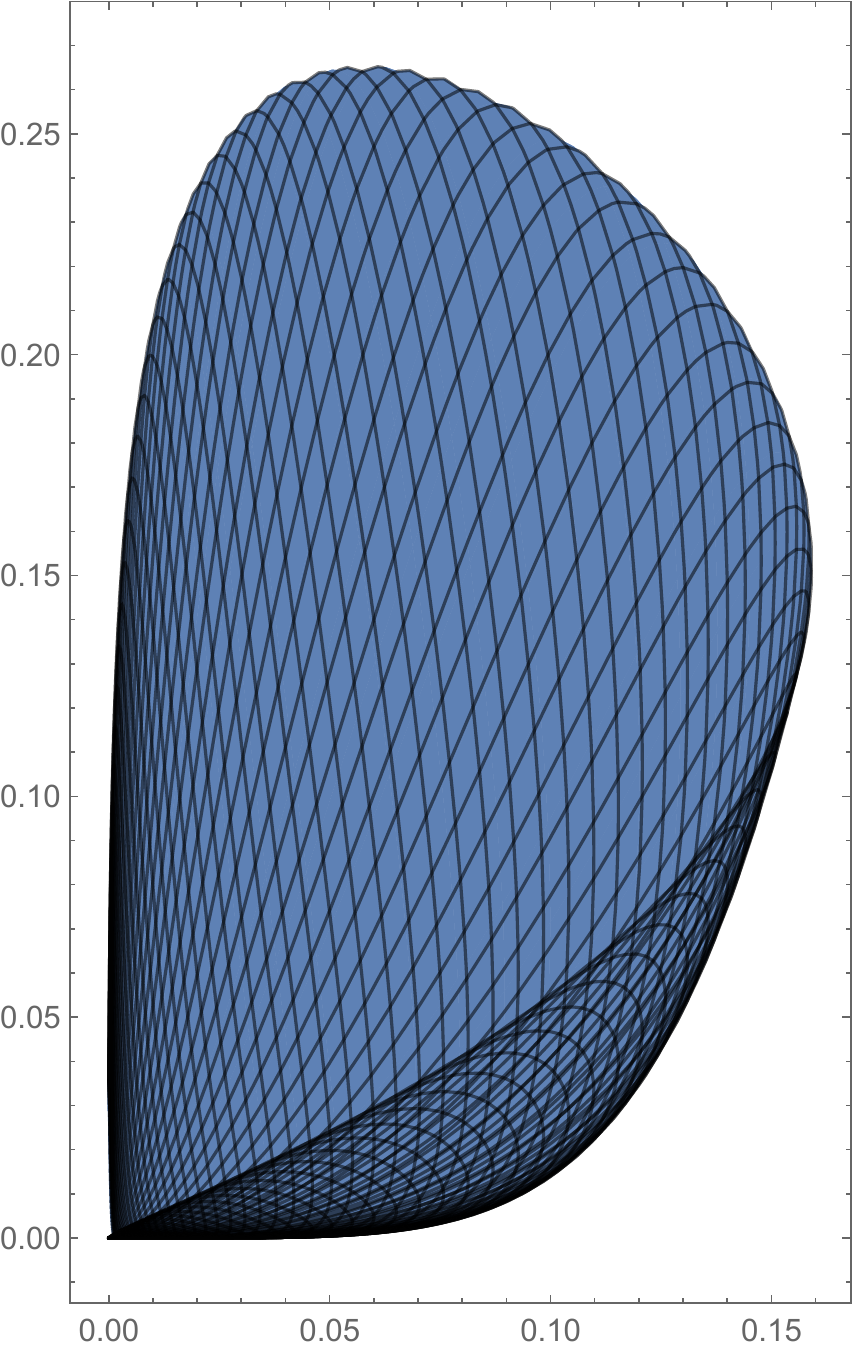}
\end{center}
\caption{Type 2: The left figure shows the contours of the densities $f_1$ (blue) and $f_2$ (red) which are codirectional. Here the black curves are the singular set of $F=(f_1,f_2)\colon\R^2\to\R^2$. The right figure shows the image of the mapping $F\colon\R^2\to\R^2$.}
\label{fig:type2}
\end{figure}

\begin{exam}
Let
$\bmu_1=(0,0)$, 
$\Sigma_1=\left(\begin{array}{cc}1&0\\0&1\\\end{array}\right)$,
$\bmu_2=(1,1)$, 
$\Sigma_2=\left(\begin{array}{cc}1&0.8\\0.8&1\\\end{array}\right)$.
Then the pair $(f_1, f_2)$ is of Type 2. Figure \ref{fig:type2}  shows the contours of $f_1,f_2$ and the image of the mapping $F\colon\R^2\to\R^2$. In particular, $S(F)$ is two intersecting lines, and $F$ is locally $\A$-equivalent to the normal form $(x,xy^2+y^4)$ at the node point.
\end{exam}

\subsection{Type 3}
\begin{figure}[htbp]
\begin{center}
\includegraphics[clip,width=6.0cm]{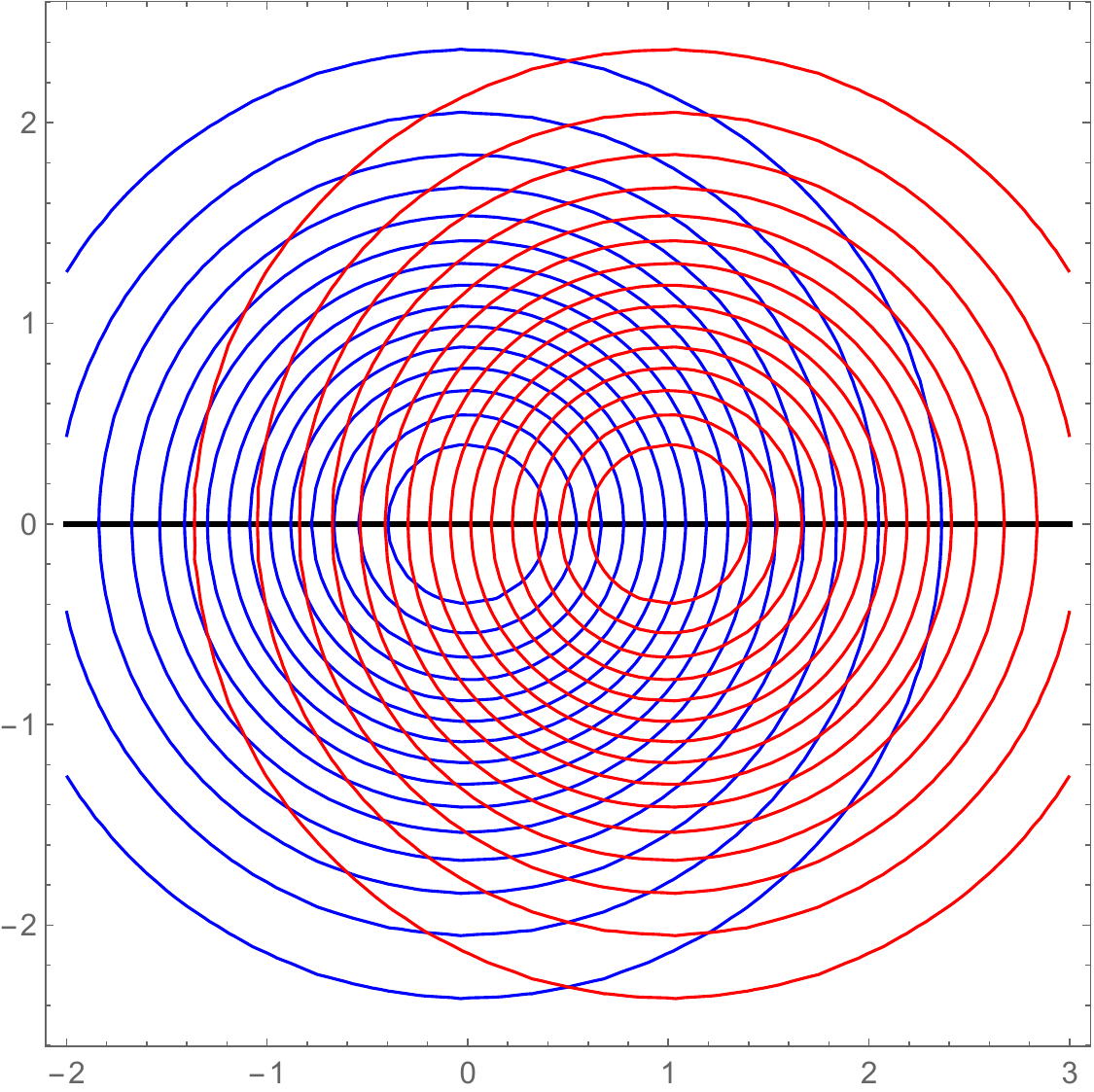}\;\;
\includegraphics[clip,width=6.0cm]{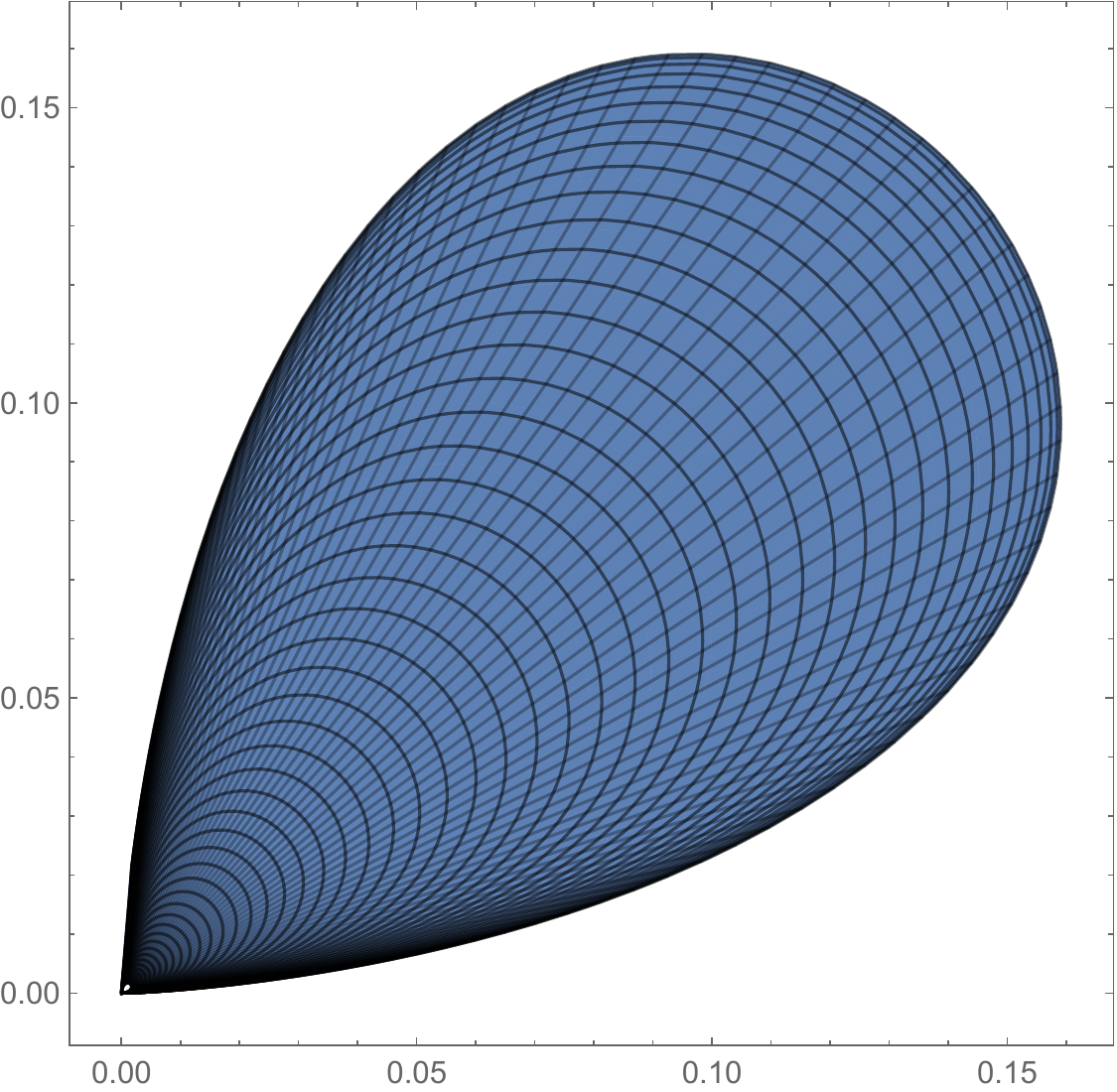}
\end{center}
\caption{Type 3: The left figure shows the contours of the densities $f_1$ (blue) and $f_2$ (red) with the variances $\Sigma_1, \Sigma_2$ being proportional. Here the black curve is the singular set of $F=(f_1,f_2)\colon\R^2\to\R^2$. The right figure shows the image of the mapping $F\colon\R^2\to\R^2$.}
\label{fig:type3}
\end{figure}

\begin{exam}
Let
$\bmu_1=(0,0)$, 
$\Sigma_1=\left(\begin{array}{cc}1&0\\0&1\\\end{array}\right)$,
$\bmu_2=(1,0)$, 
$\Sigma_2=\left(\begin{array}{cc}1&0\\0&1\\\end{array}\right)$.
Then the pair $(f_1, f_2)$ is of Type 3. Figure \ref{fig:type3}  shows the contours of $f_1,f_2$ and the image of the mapping $F\colon\R^2\to\R^2$. In particular, $S(F)$ is a line, and any singularity of $F$ is a fold.
\end{exam}

\bigskip
\noindent
{\bf Acknowledgement}. 
This work is partially supported by JSPS KAKENHI Grant Number JP 20K14312, JP 21K18312 and 24K22308.
The authors wold like to thank Toshizumi Fukui, Shunsuke Ichiki and Satoshi Kuriki for helpful comments and discussions.

\noindent
Yutaro Kabata,\\
School of Information and Data Sciences, 
Nagasaki University, 
Nagasaki 852-8521, Japan.\\
E-mail address: kabata@nagasaki-u.ac.jp
\\
\\
Hirotaka Matsumoto,\\ 
School of Information and Data Sciences,
Nagasaki University, 
Nagasaki 852-8521, Japan.
\\
\\
Seiichi Uchida,\\
Department of Advanced Information Technology, Kyushu University, 744 Motooka, Nishi-Ku, Fukuoka 819-0395, Japan
\\
\\
Masao Ueki,\\
School of Information and Data Sciences,
Nagasaki University, 
Nagasaki 852-8521, Japan.

\end{document}